\documentstyle[amsmath,amssymb,amscd,graphicx,psfrag,epsf,amsthm]{article}

\newtheorem{thm}{Theorem}

\newtheorem{lem}[thm]{Lemma}
\newtheorem{cor}[thm]{Corollary}

\newtheorem{prop}[thm]{Proposition}

\newtheorem{conj}[thm]{Conjecture}
   
\theoremstyle{definition}

\newtheorem{say}[thm]{}
\newtheorem{exmp}[thm]{Example}

\newtheorem{rem}[thm]{Remark}          
\newtheorem{ack}{Acknowledgments}

\newtheorem{defn-thm}[thm]{Definition--Theorem}  

\theoremstyle{remark}
\newtheorem{claim}[thm]{Claim}

\renewcommand{\c}[0]{{\mathbb C}}  

\renewcommand{\o}[0]{{\mathcal O}} 

\renewcommand{\r}[0]{{\mathbb R}} 

\renewcommand{\a}[0]{{\mathbb A}}

\newcommand{\p}[0]{{\mathbb P}}
\newcommand{\f}[0]{{\mathbb F}}
\newcommand{\q}[0]{{\mathbb Q}}
\newcommand{\map}[0]{\dasharrow}
\newcommand{\qtq}[1]{\quad\mbox{#1}\quad}
\newcommand{\spec}[0]{\operatorname{Spec}}

\newcommand{\gal}[0]{\operatorname{Gal}}

\newcommand{\im}[0]{\operatorname{im}}

\newcommand{\Hom}[0]{\operatorname{Hom}}

\newcommand{\onto}[0]{\twoheadrightarrow}

\def\into{\DOTSB\lhook\joinrel\rightarrow}

\begin{document}
\bibliographystyle{amsplain}

\title{Rationally connected varieties and fundamental groups}
\author{J\'anos Koll\'ar}

\maketitle

The Lefschetz hyperplane theorem says that if
$X\subset \p^N$ is a smooth projective variety over $\c$ and
$C\subset X$ is a smooth curve which is a complete intersection
of hyperplanes with $X$ then
$$
\pi_1^{top}(C)\to \pi_1^{top}(X)\qtq{is surjective,}
$$
where $\pi_1^{top}$ is the topological fundamental group.
Later  a quasiprojective version of this
result was also established (see, for instance, \cite{gmp88}
for a discussion and further generalizations).
 This says  that if $X^0\subset X$ is open 
and $C\subset X$ as above is 
{\it sufficiently general},
then
$$
\pi_1^{top}(C\cap X^0)\to \pi_1^{top}(X^0)\qtq{is surjective.}
$$

On a rationally connected variety one  would like
to use rational curves to obtain a similar result.
Complete intersection curves are essentially never 
rational. (For instance, if $X\subset \p^n$ is a hypersurface
then a general complete intersection with hyperplanes is
a rational curve iff $X$ is a hyperplane or a quadric.)
Therefore we have to proceed in a quite different way.

Let $X$ be a smooth, projective, rationally connected  variety over $\c$
and $X^0\subset X$ an open set whose complement is a normal crossing
divisor $\sum D_i$. Let $C\subset X$ be  a smooth rational curve
which intersects  $\sum D_i$ transversally everywhere. 
It is then easy to see that the
{\em normal} subgroup
 of $\pi_1^{top}(X^0)$ generated by the image of $\pi_1^{top}(C\cap X^0)$
is $\pi_1^{top}(X^0)$ itself. 
One can also achieve that the image of $\pi_1^{top}(C\cap X^0)$
has finite index in $\pi_1^{top}(X^0)$. The map
$\pi_1^{top}(C\cap X^0)\to \pi_1^{top}(X^0)$ is, however, not always  onto. 
Some examples are given in \cite[4.5]{koll-rcfg}.

$\pi_1^{top}(X^0)$ is typically an infinite group and in many questions
the difference between it and its finite index subgroups is minor.
However, in the arithmetic applications, for instance in 
(\ref{tors.gen.cor}) and
in \cite{ko-sz} it is crucial to prove surjectivity.

In order to work in arbitrary characteristic,
we should use  the algebraic fundamental group, which we denote by $\pi_1$.
(A very good introduction to algebraic fundamental groups
is \cite{mezard}. A more thorough treatment is given in
\cite{murre} and the ultimate reference is
 \cite{sga1}.)
All the crucial points of the arguments
in this paper  can be seen by concentrating on
the case of complex varieties and the topological fundamental group.

We  care only about finite quotients of the
algebraic fundamental group. Classically, 
quotients of the fundamental group correspond to
covering spaces, 
and similarly, finite \'etale covers correspond to
 finite index  subgroups
of the algebraic fundamental group. 
(Strictly speaking, the correspondence is with
finite index  open subgroups, but
the relevant topology plays no essential role in what follows.)
 A consequence of this
is that a surjectivity between fundamental groups
can be checked on finite quotients. 
We adopt this as a definition:

\begin{defn-thm} \label{basic.onto.lem} (cf.\ \cite[p.94]{murre})
Let $k$ be a field and $p:Y\to X$  a morphism of normal, 
geometrically connected $k$-schemes.
Then
$\pi_1(Y)\to \pi_1(X)$ is {\it surjective} if
for every finite \'etale cover $X'\to X$ with $X'$ connected,
the pull back $Y\times_XX'$ is also connected.

Similarly, the image of
$\pi_1(Y)\to \pi_1(X)$ is said to have  index $\leq m$ if
for every finite \'etale cover $X'\to X$ with $X'$ connected,
the pull back $Y\times_XX'$ has $\leq m$ connected components.

There are infinite groups without finite
quotients. Thus there are compact manifolds which are
not simply connected but have no finite sheeted covering
spaces.  It is not known if this also occurs
for smooth projective varieties over $\c$. In any case,
the algebraic  fundamental group detects
only the finite quotients of the topological fundamental group
in characteristic zero.
\end{defn-thm}

\begin{rem}[Base points]
As in topology, the fundamental group can be defined
only with a fixed base point.
Fundamental groups of the same scheme with different
base points are isomorphic, but the isomorphism is well
defined only up to inner automorphisms.
For some questions, for instance for surjectivity
of maps, this does not matter. In such cases we 
occasionally omit the base point from the notation.

It is, however, important to keep in mind  that ignoring the base
point can easily lead to wrong conclusions.

In algebraic geometry the base point is 
 a {\it geometric point} of $X$. That is, a  morphism
$\bar x\to X$ where $\bar x$ is the spectrum of an algebraically closed field.
If $x\in X$ is a point and $K$ is any algebraically closed field containing 
the residue field $k(x)$, then we get a geometric point $\bar x:=\spec K$
with a map $\bar x\to X$.
 The resulting fundamental group $\pi_1(X,\bar x)$ does not depend on
the choice of $K$ and we usually denote it by $\pi_1(X,x)$. 

As in topology, we have the equivalence
$$
\begin{array}{rcl}
\left\{
\begin{array}{c}
\mbox{finite index open}\\
\mbox{subgroups of $\pi_1(X,x)$}
\end{array}
\right\}
&\Leftrightarrow&
\left\{
\begin{array}{c}
\mbox{connected finite \'etale}\\
\mbox{covers of $X$ plus a}\\ 
\mbox{lifting of the base point}
\end{array}
\right\}
\end{array}
$$

\end{rem}

\section{Statements of the main results}

Separably rationally connected  (or SRC) 
 varieties were 
introduced in \cite{kmm2}. See   \cite[34]{arko}
in this volume
for the definition and basic properties.
(Note in particular that an SRC variety is always 
assumed to be geometrically connected.)
A weaker version of 
our main Theorem (\ref{ratlefx.thm})  was proved
in \cite{koll-rcfg}. 
The key improvement is that the present proof, besides
being simpler, does not need the resolution of singularities,
thereby extending the range of applications to
positive characteristic. It is also nicer  that
the choice of the curve $C$ is largely independent of
the open set $X^0$.

We formulate 
the main result in two versions,  with
fixed base point and with variable base point.
 Both are quite useful and it seems
very cumbersome to have a unified statement.

\begin{thm} \label{ratlefx.thm}
Let $X$ be a smooth, projective, SRC variety over a field $k$
and $x\in X(k)$.
 Then
 there is a dominant family of rational curves through $x$
(defined over $k$)
$$
F:W\times \p^1\to X,  \quad F(W\times (0{:}1))=\{x\}
$$
with the following properties:
\begin{enumerate}
\item $W$ is geometrically irreducible, smooth
and open in $\Hom(\p^1, X, (0{:}1)\mapsto x)$ and
$F: W\times (1{:}0)\to X$ is  smooth.
\item There is a $k$-compactification $\bar W\supset W$
such that $\bar W$ has a smooth $k$-point.
\item $F_w^*T_X$ is ample for every $w\in W$.
\item Let $K\supset k$ be an algebraic closure and $X^0_K\subset X_K$
any open $K$-subvariety containing $x$.  Then 
 for every $w\in W(K)$ the induced map on the geometric fundamental groups
$$
\pi_1(F^{-1}(X^0_K), (w,(0{:}1)))\onto
\pi_1(X^0_K,x)
\qtq{is surjective.}
$$
\end{enumerate}
\end{thm}

\begin{rem} \label{main.thm.rem}
From the proof we see that we have considerable freedom in constructing
$W$ and $F$. Thus we can construct them with additional properties
that are useful in some applications.
\begin{enumerate}\setcounter{enumi}{4}
\item If $\dim X\geq 3$ then
we can also assume that $F$ is an embedding on every $\p^1_w$.
If $\dim X =2$ then 
we can also assume that $F$ is an immersion  on every $\p^1_w$.
\item Given an integer $d$, we can also assume that
 $F_w^*T_X(-d)$ is ample for every $w\in W$.
\item Fix an embedding $X\subset \p^N$. Then the degrees
of the curves $F(\p^1_w)$ can be bounded from above in terms
of the degree and dimension of $X$  (and the above integer $d$).
\end{enumerate}
We point out at the end of the proof 
(\ref{main.thm.rem.pf}) how to achieve these
additional properties.

These comments also apply to (\ref{ratlef.thm}).

\end{rem}

 Based on topology, one would expect that
for $w\in W(K)$ general, the induced map
$$
\pi_1(\p_w^1\cap F^{-1}(X^0_K),(0{:}1))\to
\pi_1(X^0_K,x)
$$
is surjective. This is indeed so in characteristic zero,
but in characteristic $p$ the \'etale covers of $p$-power degree
behave unexpectedly (see, for instance, example  (\ref{chrp.strange.exmp})). 
The following weaker result, however, still holds, and this is
enough for most applications.

\begin{cor}  \label{ratlefx.cor}
Notation and assumptions as in (\ref{ratlefx.thm}).
Let $\pi_1(X^0_K,x)\onto G$ be a finite quotient
(whose kernel is open).
Then there is
a dense open subset
$W^0_G\subset W$ such that for every $w\in W^0_G(K)$ the composed map
$$
\pi_1(\p_w^1\cap F^{-1}(X^0_K),(0{:}1))\to
\pi_1(X^0_K,x)\to G 
\qtq{is surjective.}
$$
\end{cor}

The base point free version of (\ref{ratlefx.thm})
is the following:

\begin{thm} \label{ratlef.thm}
Let $X$ be a smooth, projective SRC variety over a field $k$.
 Then
 there is a  dominant family of rational curves
(defined over $k$)
$$
F:W\times \p^1\to X
$$
with the following properties:
\begin{enumerate}
\item $W$ is geometrically irreducible and smooth.
\item $F$ is smooth and $F_w^*T_X$ is ample for every $w\in W$.
\item Let $K\supset k$ be an algebraic closure and $X^0_K\subset X_K$
any open $K$-subvariety.  Then 
 for every $z\in  F^{-1}X^0_K(K)$
 the induced map on the geometric fundamental groups
$$
\pi_1(F^{-1}X^0_K, z)\onto
\pi_1(X^0_K,F(z))
\qtq{is surjective.}
$$
\end{enumerate}
\end{thm}

\begin{cor}\label{ratlef.cor2}  
Notation and assumptions as in (\ref{ratlef.thm}).
Let $Y_K$ be an irreducible $K$-variety and
$h:Y_K\to X_K$ a dominant morphism.

Then there is a dense open set $W^0_Y\subset W_K$
such that for every  $w\in W^0_Y(K)$  
the fiber product  $Y_K\times_{X_K} \p^1_w$
 is  irreducible.
\end{cor}

\section{Fundamental groups in families}

In this section we gather some 
fundamental group lemmas. The main result is
(\ref{gen.semicont.cor}) comparing the fundamental groups of
different fibers of a morphism.

\begin{lem}\label{open.onto.lem}
 Let $X$ be a normal variety and $X^0\subset X$ an open
subvariety. Then for every geometric point $x\to
 X^0$
the induced map $\pi_1(X^0,x)\onto \pi_1(X,x)$ is surjective.
\end{lem}

Proof. Let $p:Y\to X$ be a finite \'etale cover with $Y$ connected.
Since $Y$ is normal, it is also 
irreducible.
Then  $p^{-1}(X^0)\subset Y$ is a dense open subscheme,
hence also irreducible and connected.\qed

\begin{lem} \label{npstein.lem}
Let $p:Y\to X$ be a dominant morphism of
normal varieties. 
Then there is an  open subset $X^0\subset X$ and a factorization
$$
p^{-1}(X^0)\stackrel{q}{\to}  Z^0\stackrel{r}{\to} X^0
$$
such that $q$ has only geometrically irreducible fibers
and $r$ is finite and \'etale. $Z^0$ is also normal
and if $Y$  is irreducible then so is $Z^0$.
\end{lem}

Proof. Let $\bar Y\supset Y$ be a normal variety such that
$p$ extends to  a proper morphism $\bar p:\bar Y\to X$.
Let $\bar Y\to W\to X$ be the Stein factorization of $\bar p$.
$W\to X$ is finite and dominant, hence $W$ is the normalization
of $X$ in the function field $k(W)$. 
Let $F\subset k(W)$ denote the separable closure of $k(X)$ in $k(W)$.
Let $Z$ be the normalization
of $X$ in the function field $F$. Then $k(W)/k(Z)$ is purely inseparable and
$k(Z)/k(X)$ is separable. Thus
$Z\to X$ is  \'etale over an open set $X^1\subset X$.

Next we claim that
there is an open subset $W^1\subset W$ such that the fibers
of $\bar Y\to W$ are geometrically irreducible over $W^1$. This is
a little tricky in positive characteristic.
Set $K=k(W)$. 

The generic fiber  $\bar Y_K$ is normal
and $H^0(\bar Y_K,\o)=K$ by the definition of Stein factorization.
Unfortunately this does not imply that the geometric generic fiber
$\bar Y_{\bar k}$ is also normal.
Nonetheless,  by 
 (\ref{normal.extend.lem}) 
the geometric generic fiber  $\bar Y_{\bar K}$ is irreducible.

Thus there is an open subset $W^1\subset W$ such that the fibers
of $\bar Y\to W$ are geometrically irreducible
 over $W^1$ (cf.\ \cite[IV.12]{ega}). Let $X^2\subset X$ be an open subset
 whose preimage in $W$ is contained in $W^1$.
Taking $X^0:=X^1\cap X^2$ we are done.\qed

\begin{lem} \label{normal.extend.lem}
 Let $Z_K$ be a normal proper
variety over a field $K$ and assume that $H^0(Z,\o_Z)=K$.

Then  $Z$ is geometrically irreducible, that is, 
$Z_{\bar K}$ is irreducible.
\end{lem}

Proof.  
If $L/K$ is a separable field extension then
$Z_L:=Z_K\times_K\spec L$ is also normal
(see, for instance \cite[21.E]{mats80}). If we are in characteristic zero
then $\bar K/K$ is separable, thus
  $Z_{\bar K}$ is also normal.
Cohomologies of sheaves change by tensoring under
field extensions, thus $H^0(Z_{\bar K},\o_{Z_{\bar K}})=\bar K$,
and so $Z_{\bar K}$ is connected. A normal and connected scheme is irreducible.

In positive characteristic we have a
 problem  with inseparable extensions of $K$.

Let
 $L/K$ be  a finite degree  purely inseparable extension. 
Let $Z_L^n$ denote the normalization of $Z_L$.
Then $Z_L^n$ is also the normalization of $Z_K$
in the function field $L(Z)$. $L(Z)/K(Z)$
is purely inseparable, thus there is a $p$-power $q$
such that the $q$-power map
$F_q:g\mapsto g^q$ maps $L(Z)$ into $K(Z)$.
This implies that, up to  a power of the Frobenius map,
$Z_L^n\to Z_K$ can be inverted. In particular,
$Z_L^n\to Z_K$ is a {\em homeomorphism}.

Every finite extension of $K$ can be written as  a separable extension
followed by an inseparable one.
Putting these two steps together we get that if
$Y_K$ is normal and
geometrically connected  then  it is
geometrically irreducible. \qed

\begin{lem}\label{dom.find.lem}
 Let $p:Y\to X$ be a dominant morphism of
geometrically irreducible 
normal varieties over an algebraically closed  field $k$.
 Let $\bar y\to Y$ be a geometric point of $Y$ and $\bar x=p(\bar y)$.
Then 
$\im[\pi_1(Y,\bar y)\to \pi_1(X,\bar x)]$
is a  finite index closed subgroup of  $\pi_1(X,\bar x)$.
\end{lem}

Proof.   Let 
$
Y^0:=p^{-1}(X^0)\stackrel{q}{\to}  Z^0\stackrel{r}{\to} X^0
$
be as in (\ref{npstein.lem}).
$\pi_1(Z^0)$ has finite index closed image in $\pi_1(X^0)$,
which surjects onto $\pi_1(X)$. Thus it is
enough to prove that $\pi_1(Y^0)$ surjects onto
$\pi_1(Z^0)$. 

Let $Z^1\to Z^0$ be a connected \'etale cover.
Then $Y^0\times_{Z^0}Z^1\to Z^1$ has geometrically
irreducible fibers, hence the fiber product is also connected.
Thus  $\pi_1(Y^0)$ surjects onto
$\pi_1(Z^0)$ according to the definition  (\ref{basic.onto.lem}).\qed
\medskip

The following result was stated by \cite{campana}
over $\c$ but the proof applies in any characteristic:

\begin{cor} Let $X$ be a normal,  proper
rationally connected  variety
over an algebraically closed field $k$.
Then $\pi_1(X)$ is finite. 
\end{cor}

Proof.  
Let $x\in X$ be any point.
Since $X$ is rationally connected, there is a dominant morphism
$$
F:W\times \p^1\to X \qtq{such that} F(W\times (0{:}1))=\{x\}.
$$
By (\ref{dom.find.lem}) the image of
$\pi_1(W\times \p^1)$ has finite index in $\pi_1(X)$.
On the other hand, $\pi_1(W\times \p^1)=\pi_1(W\times(0{:}1))$
and the latter is mapped to the identity. Thus
$\pi_1(X)$ is finite.\qed
\medskip

In the smooth case, we even have simple connectedness.
In characteristic zero this 
was proved by  \cite{campana}, it
is also explained  in
\cite{deb}. The positive characteristic case follows from  a recent
result of \cite{jo-st},  a proof is written up in 
\cite[3.6]{deb-b}.

\begin{thm}\cite{kol-let} \label{SRC.sc}
 Let $X$ be a smooth projective SRC variety
over an algebraically closed field $k$ of arbitrary characteristic.
 Then
$\pi_1(X)=\{1\}$. \qed
\end{thm}

\begin{say}[Fundamental groups in families]

We also need to know how the fundamental group varies in families.
In characteristic zero the question is topological and
easy. There are some new twists in
positive characteristic. Several of these are
illustrated by the following example:

\begin{exmp} \label{chrp.strange.exmp}
Let us work over an algebraically closed field of
 characteristic $p$. Consider the surface 
$S:=(f=0)\subset \a^3$ where $f=y+z-xz^p+z^{2p}$. The derivative
$\partial f/\partial z=1$, so the projection
onto the $z=0$ plane $S\to \a^2$
is finite and \'etale. (This already shows
that $\pi_1(\a^2)$ is nontrivial.)

Consider the line $L_c:=(y-cx=0)$. The preimage of
$L_c$ in $S$ is isomorphic to the curve
$cx+z-xz^p+z^{2p}=0$. If $c=0$ then this curve
is reducible but if $c\neq 0$ then it is irreducible.
This shows that the image of
$$
\pi_1(\a^1)\cong \pi_1(L_c)\to \pi_1(\a^2)
$$
depends on $c$, and in fact the map is not surjective for $c=0$.
Since any two lines are equivalent under automorphisms of $\a^2$,
we  get that for any line $L\subset \a^2$, the
natural map
$$
\pi_1(\a^1)\cong \pi_1(L)\to \pi_1(\a^2)
$$
is not surjective.
It is not hard to see, however, that the following slightly
weaker version holds:

Let $\pi_1(\a^2)\to G$ be a finite quotient (with open kernel).
Then the composite
$$
\pi_1(\a^1)\cong \pi_1(L)\to \pi_1(\a^2)\to G
$$
is surjective for a general line $L$, the notion of
general depending on $G$. 
\end{exmp}

The latter observation holds very generally:

\begin{prop} \label{gen.semicont.prop} Let $k$ be an 
algebraically closed field.
Let $f:Y\to X$ be a morphism of  irreducible and normal $k$-varieties.
Let $S\subset Y$ be an irreducible sub-variety
such that $S\to X$ is dominant.
For $s\in S(k)$ let $Y_s$ denote the fiber of $f$ over $f(s)$.
Let $Z$ be a variety, $z\in Z(k)$  a point and $\pi_1(Z,z)\to G$
a finite quotient (with open kernel).
 Let $h:Y\to Z$ be a morphism
such that $h(S)=\{z\}\in Z$. 

Then there is a dense open subset $S^0\subset S$
such that for every $s\in S^0(k)$ the two maps
$$
\pi_1(Y_s,s)\to \pi_1(Z,z)\to G
\qtq{and}
\pi_1(Y,s)\to \pi_1(Z,z)\to G
$$
have the same image.
\end{prop}

Proof. We are allowed to change base from $X$ to $S$
and then replace $S$ by an open subset. Thus we can assume that
$S\cong X$ is a section.

Let $Z'\to Z$ be the cover corresponding 
to the image of $\pi_1(Y,s)\to \pi_1(Z,z)\to G$.
We need to prove that 
the number of irreducible components of $Z'\times_ZY$
is the same as the number of irreducible components of
$Z'\times_ZY_s$  for $s\in S^0$.
Let $W\subset Z'\times_ZY$ be an irreducible component. Then $W$ is
an irreducible \'etale cover of $Y$ which is trivial
over $S$. Let $g:W\to X$ be the composite.

As in (\ref{npstein.lem}) there are  open sets $X^0\subset X$, $W^0\subset W$
and a factorization  $W^0\to V^0\to X^0$ 
where $W^0\to V^0$ has geometrically irreducible fibers
and $V^0\to X^0$ is finite and \'etale. Since $S\subset W$
is a section of $g$, we have $X^0\cong S\cap W^0\to V^0\to X^0$,
thus $V^0\cong X^0$. Therefore $W^0\to X^0$ 
has geometrically irreducible fibers.
\qed
\medskip

As a corollary we obtain a lower semicontinuity statement
for fundamental groups of fibers.
We do not compare the groups themselves, just their images
in finite groups whenever this can be done sensibly.

\begin{cor}\label{gen.semicont.cor} Let $k$ be an 
algebraically closed field.
Let $f:Y\to X$ be a morphism of irreducible and 
 normal $k$-varieties with a section
 $S\subset Y$.
For $s\in S$ let $Y_s$ denote the fiber of $f$ over $f(s)$.
Let $Z$ be a $k$-variety, $z\in Z(k)$  a point 
 and  $h:Y\to Z$  a morphism
such that $h(S)=\{z\}\in Z$.
\begin{enumerate}
\item For any $s,s'\in S$,
$$
\im[\pi_1(Y_s,s)\to \pi_1(Z,z)]
\subset
\im[\pi_1(Y,s')\to \pi_1(Z,z)].
$$
\item  Let $S^0\subset S$ be any  dense open subset. Then
$$
\im[\pi_1(Y,s')\to \pi_1(Z,z)]=
\langle
\im[\pi_1(Y_s,s)\to \pi_1(Z,z)]: \forall s\in S^0
\rangle.
$$
\end{enumerate}
\end{cor}

Proof. The maps $Y_s\to Y\to Z$ give the inclusion
$$
\im[\pi_1(Y_s,s)\to \pi_1(Z,z)]
\subset
\im[\pi_1(Y,s)\to \pi_1(Z,z)].
$$
Since $S$ is irreducible and contracted to $z$,
$$
\im[\pi_1(Y,s)\to \pi_1(Z,z)]
=
\im[\pi_1(Y,s')\to \pi_1(Z,z)]
$$
for any $s,s'\in S$, proving (1). 

This implies the containment $\supset$ in (2), and we only need
the reverse inclusion.
By the correspondence between the fundamental group and
finite \'etale covers we need to prove that
if $Z'\to Z$ is a finite  \'etale cover 
such that $Y_s\times_ZZ'$ is reducible for every $s\in S^0$
then $Y\times_ZZ'$ is also reducible. This is
so by (\ref{gen.semicont.prop}).\qed

\end{say}

\begin{rem}\label{fg.of.combs}
Several times we will apply (\ref{gen.semicont.cor})
to smoothings of combs, as defined in \cite[43]{arko}.
Let $C=C_0\cup\dots\cup  C_n$ be a pointed comb with handle $C_0\ni p$.
Let $c_i\in C_i$ be the nodes and  $C^*\subset C$  an open
subset containing $C_0$. 
The injections $C_i\into C$ induce maps 
$\pi_1(C^*\cap C_i, c_i)\to \pi_1(C^*,c_i)$. 
(These are in fact injections, but this is not crucial to us.)
Thus we obtain maps $\pi_1(C^*\cap C_i, c_i)\to \pi_1(C^*,p)$
which are well defined modulo conjugation by an element
of $\pi_1(C_0,p)$.

In particular, if $F:C\to X$ is 
assembled from the maps $f_i:C_i\to X$
and $X^0\subset X$ is an open subset containing $x$  then
$$
\im[\pi_1(F^{-1}(X^0),p)\to \pi_1(X^0,x)]
$$
contains  all the images
$$
\im[\pi_1(f_i^{-1}(X^0),c_i)\to \pi_1(X^0,x)]
\qtq{for $i=1,\dots, n$.}
$$
At least when the base field is algebraically closed, 
the conjugation by elements
of $\pi_1(C_0,p)$ does not matter since 
$F(C_0)=\{x\}$ so $\pi_1(C_0,p)$ is mapped to the identity.

Note that 
a similar assertion would not hold
without the assumption $F(C_0)=x$.
\end{rem}

\section{Proofs of the theorems}

First we reduce  (\ref{ratlefx.thm}) 
to the a priori much weaker
(\ref{weakratlef}). Then in the second part we
prove (\ref{weakratlef}).
Finally we explain that  (\ref{ratlefx.thm}) 
implies the other statements.
We repeatedly use several properties of
stable curves discussed in \cite[Sec.8]{arko}.

\begin{prop} \label{weakratlef}
Let $X$ be a smooth, projective, SRC variety over an
algebraically closed  field $k$
and $x\in X(k)$.
 Then
 there is a 
family of free genus zero stable maps with base point $x$ 
$$
(C/S,F:C\to X\times S, \sigma:S\to C)
$$ 
parametrized by a scheme of finite type $S$
such that
for every open subset $X^0\subset X$ with $x\in X^0$
the images
$$
\im[\pi_1(F_s^{-1}(X^0), \sigma_s(s))\to
\pi_1(X^0,x)]\qtq{for all $s\in S$}
$$
(topologically)
generate $\pi_1(X^0,x)$.
\end{prop}

\begin{say}[Proof of (\ref{weakratlef})$\Rightarrow$(\ref{ratlefx.thm})]
\label{pf.of.23->4}

Assume first that $k$ is algebraically closed.
We start with the family
$(C/S,F:C\to X\times S, \sigma:S\to C)$ in (\ref{weakratlef})
and using it we repeatedly construct other families
until we end up with one as required for (\ref{ratlefx.thm}).

$(C/S,F:C\to X\times S, \sigma:S\to C)$ 
 gives the moduli  morphism
$S\to \bar M_0(X,\spec k\to x)$ (cf. \cite[41]{arko});
 let $S^*$ be its image (as a constructibe set).

Pick any $s\in S(k)$ and let 
$$
\begin{array}{ccccc}
C_s &\subset&  C_{U_s} & \stackrel{\Phi_{U_s}}{\to}  X\\
\downarrow &&  \downarrow  &&\\
 s& \in & U_s &&
\end{array}$$
be a  versal  family as in \cite[42.5]{arko}.

To each $U_s$ there corresponds a natural morphism
$U_s\to \bar M_0(X,\spec k\to x)$ with open image.
Finitely many of these cover $S^*$. From now on
we work with  these finitely many families
$$
(C_{U_s}/U_s,\Phi_{U_s}:C_{U_s}\to X\times U_s,\sigma_{U_s}:U_s\to C_{U_s}).
$$

From (\ref{gen.semicont.cor}) we see that 
the conclusion of (\ref{weakratlef}) remains true if we
replace $\cup U_s$ by a dense open subset.
By \cite[42.3]{arko}
 the points corresponding to irreducible curves are dense
in each $U_s$, thus we may replace each $U_s$ by the
open subset of irreducible and free curves.

We now have finitely many (say $N$) families 
of stable and free maps
$$
\Phi_j:U_j\times \p^1\to X\qtq{with}  \Phi_j(U_j\times(0{:}1))=x,
$$
where each $U_j$ is irreducible. 
We are free to add one more family $\Phi_{N+1}:U_{N+1}\times \p^1\to X$
consisting of very free curves through $x$.
In order to create a single family, we use a trick introduced in \cite{rcloc}.

Pick points $u_j\in U_j$ and, as in \cite[43]{arko},
 assemble a pointed comb
out of the maps  $\Phi_{j,u_j}:\p^1\to X$.
Since  $N+1\geq 2$,  these pointed combs are stable genus zero curves
over $X$.
All these combs are parametrized by an irreducible
variety  (the product of the $U_j$ times the space of
$N+1$ points on $\p^1$). Thus all these combs
are in the same irreducible component 
${\mathcal U}\subset \bar M_0(X,\spec k\to x)$.

As before, using 
(\ref{fg.of.combs}) we can first replace
the families $\Phi_j:U_j\times \p^1\to X$
with
the irreducible component ${\mathcal U}$ of
the  family of stable curves 
containing all of the above combs 
$$
f_C:C=C_0\cup\cdots \cup  C_{N+1}\to X,
$$
and then by (\ref{gen.semicont.cor}) with  the open 
subfamily  ${\mathcal U}^0\subset {\mathcal U}$
of
irreducible  and free curves. The moduli map
${\mathcal U}^0\to \Hom(\p^1, X, (0{:}1)\mapsto x)$
is dominant onto one of the irreducible components.
Let $V\subset \Hom(\p^1, X, (0{:}1)\mapsto x)$
be an open set contained in the image of ${\mathcal U}^0$.

Thus we obtain a single irreducible family of stable maps
$$
\Phi:V\times \p^1\to X\qtq{with}  \Phi(V\times(0{:}1))=x.
$$
We claim that this
 satisfies all the requirements of
(\ref{ratlefx.thm}) for $k$ algebraically closed.

Indeed, 
 $V$ is geometrically irreducible, smooth
and open in $\Hom(\p^1, X, (0{:}1)\mapsto x)$
by construction.
By \cite[II.3.4]{Ko96}, 
the tangent map of $F:(V\times (1{:}0))\to X$ at $(v,(1{:}0))$
is given by the restriction
$$
H^0(\p^1, F_v^*T_X(-(0{:}1)))\to  F_v^*T_X|_{(1{:}0)}\cong T_{X,F(v,(1{:}0))}
$$
which is surjective if $H^1(\p^1, F_v^*T_X(-2))=0$. We prove the latter below,
thus $F:(V\times (1{:}0))\to X$ is smooth, at least after
passing to a suitable open subset of $V$.

Condition (\ref{ratlefx.thm}.2) is
not interesting when $k$ is algebraically closed.

Ampleness of a vector bundle $E$ on $\p^1$ is equivalent
to $H^1(\p^1,E(-2))=0$. By upper semicontinuity,
$H^1(\p^1, F_v^*T_X(-2))=0$ for general $v\in V$
if $H^1(C,f_C^*T_X(-2p))=0$ for one of the combs
$f_C:C\to X$ and for some smooth point
$p\in C$. Pick $p\in C_{N+1}$. By our choice of the $(N+1)$st family,
$F^*T_X$ restricted to $C_{N+1}$ is ample
and all the other restrictions  to $C_i$ are semipositive.
This easily implies that $H^1(C,f_C^*T_X(-2p))=0$,
cf. \cite[18]{arko}. This proves (\ref{ratlefx.thm}.3)
and also (\ref{ratlefx.thm}.1).

Applying (\ref{gen.semicont.cor}) several times, we get
that the closed subgroup of $\pi_1(X^0,x)$ generated by the images of
$\pi_1(\p^1_v\cap F^{-1}(X^0),(0{:}1)): v\in  V$
is $\pi_1(X^0,x)$ itself. On the other hand,
the image of 
$\pi_1(F^{-1}(X^0_K), (v,(0{:}1)))\to
\pi_1(X^0_K,x)$ is closed by (\ref{dom.find.lem}),
thus $\pi_1(F^{-1}(X^0_K), (v,(0{:}1)))\to
\pi_1(X^0_K,x)$ is surjective, giving (\ref{ratlefx.thm}.4).

Assume next that $k$ is not algebraically closed,
and let $\bar k$ be its algebraic closure.
The family constructed above is already defined
over a finite Galois extensions $k'/k$ with Galois group $G$.
(This is automatic in characteristc 0. In general,
note that free curves give smooth points of $\Hom(\p^1,X)$,
and the smooth geometric components of a $k$-scheme are defined
over a separable extension.)
We may also assume that $V$ has  a 
distinguished
$k'$-point $v\in V(k')$. By \cite[II.3.14.2]{Ko96}
we may also assume that the corresponding
map $\{v\}\times \p^1\to X$ is an immersion.

For every $\rho\in G$ let 
$$
\Phi^{\rho}:V^{\rho}\times \p^1\to X
$$
 be the conjugate family obtained by applying $\rho$
to all the defining equations.

We now proceed as above.
Pick arbitrary points $v_{\rho}\in V^{\rho}(\bar k)$
 and assemble a pointed comb
out of the maps  $\Phi^{\rho}_{v_{\rho}}:\p^1\to X$.
All these combs are parametrized by an irreducible
variety  
which is now defined over $k$.
(Over $k'$ this parameter space 
is the product of the $V^{\rho}$ times the space of
$N$ points on $\p^1$.
Over $k$, it is fibered over the Weil restriction ${\frak R}_{k'/k}V$
(see, for instance \cite[7.6]{blr})
with various forms  of the space of
$|G|$ points on $\p^1$ as fibers.)
As before, all these combs
are in the same geometrically irreducible component 
of ${\mathcal W}\subset \bar M_0(X,\spec k\to x)$.

By (\ref{gen.semicont.cor}) we can replace 
${\mathcal W}$
with  the subfamily $W$ of
irreducible  and free curves.

In order to check the required properties, 
the only new aspect is (\ref{ratlefx.thm}.2).
The comb $C$ assembled from the conjugates of $\{v\}\times \p^1$
gives a  smooth $k$-point of ${\bar{\mathcal W}}$.
Indeed, all the conjugates of  $\{v\}\times \p^1\to X$ are
immersions, thus $C$ has no automorphisms
commuting with $C\to X$ and fixing the base point.
Therefore $C$ gives a smooth point of ${\bar{\mathcal W}}$
by \cite[2.iii]{FuPa}.
\qed
\end{say}

\begin{say}[Idea of the proof of (\ref{weakratlef})]

Take all possible free, stable, genus 0 curves
$g:(c\in C)\to (x\in X)$ and 
let $H\subset \pi_1(X^0,x)$ be  the subgroup  (topologically)
generated by all the images $\pi_1(C^0,c)\to \pi_1(X^0,x)$
where $C^0:=g^{-1}(X^0)$.
We want to prove that $H=\pi_1(X^0,x)$. If not, then there is
a finite index closed subgroup $H\subset H'\subsetneq \pi_1(X^0,x)$.
The subgroup $H'$ corresponds to a degree $d$  \'etale cover
 $Y^0\to X^0$ 
together with a geometric point $y\to Y^0$ lying over $x\to X^0$.
Extend it to a finite morphism $p:Y\to X$.

Since $X$ is simply connected (\ref{SRC.sc}), $p$ is
not \'etale and so it ramifies along a divisor $D\subset Y$
(using the purity of branch loci).

Pick a general point $x'\in p(D)$ and assume that
$C\subset X$ is a rational curve passing through
$x,x'$ and transversal to $p(D)$ everywhere. We hope that
$C\to X$ cannot be lifted to $Y$.  The local picture we 
should have in mind is  
$$
Y=\c^2_{y_1,y_2}\stackrel{(y_1,y_2)\mapsto (y_1,y_2^m)}{\longrightarrow}
\c^2_{x_1,x_2}=X.
$$
Here $D=(y_2=0)$ and $p(D)=(x_2=0)$. 
Set $y':=(0,0)\in Y$.
If $C$ is the line $(x_1=0)$ through $x':=(0,0)$ 
then $p^{-1}(C)$ is the irreducible curve $(y_1=0)$ through $y'$.
Thus the inclusion $C\into X$ cannot be lifted  to $C\to Y$
in such a way that $x'\in C$ is mapped to $y'\in Y$.

This is, however, only the local picture.
The global problem is that $Y$ has many other points over $x'$
besides $y'$, and $p$
ramifies only at some of these.
The best we can say is the following.

The fiber product $C_Y:=Y\times_XC$
is a smooth curve and $C_Y\to C$
definitely ramifies at $C_Y\cap D$. Thus
$C_Y$ is {\em not} the disjoint union of  $d$ copies
of $C$. This means that there is at least one  point $y_{i_0}\in p^{-1}(x)$
such that the inclusion $C\into X$
cannot be lifted to $C\to Y$ sending $c$ to $y_{i_0}$.
For other choices $y_i\in p^{-1}(x)$ a lifting may exist.
It is very hard to tell which preimage is  which.

Choosing a preimage $y_i$
is essentially equivalent to choosing a conjugate of the
subgroup $H'$, thus we cannot move freely between the
different preimages.

We circumvent this problem as follows.

First  we choose a general smooth rational curve $B\subset X$ passing through
a general point $z\in X$ and intersecting $p(D)$ transversally everywhere.
Let $v_i\in Y$ $i=1,\dots,d$ be all the preimages of $z$.
As with $C$, there is at least one $v_{i_0}$
such that  $B\into X$ cannot be lifted
to $B\to Y$ sending $z$ to $v_{i_0}$.

Then we
 pick other rational 
 curves $A_k$ passing through $x$ and  $z$. 
We furthermore achieve (and this turns out to be  easy for a {\em general}
 point)
that the liftings of the curves $A_k$
connect $y$ with {\em all} the different points $v_i$.

This implies that at least one  of the reducible curves
$A_k\cup B$ cannot be lifted to $Y$.
(See Figure 1.)


\begin{figure}[hbtp]
\centering

\psfrag{x}{{\fontsize{20}{20}$x$}}
\psfrag{x'}{{\fontsize{20}{20}$x'$}}
\psfrag{u}{{\fontsize{20}{20}$z$}}
\psfrag{a1}{{\fontsize{20}{20}$A_1$}}
\psfrag{a2}{{\fontsize{20}{20}$A_2$}}
\psfrag{a3}{{\fontsize{20}{20}$A_3$}}
\psfrag{A1}{{\fontsize{20}{20}$A'_1$}}
\psfrag{A2}{{\fontsize{20}{20}$A'_2$}}
\psfrag{A3}{{\fontsize{20}{20}$A'_3$}}
\psfrag{b}{{\fontsize{20}{20}$B$}}
\psfrag{X}{{\fontsize{20}{20}$X$}}
\psfrag{Y}{{\fontsize{20}{20}$Y$}}
\psfrag{no}{{\fontsize{20}{20}${\rm no\ lifting}$}}
\psfrag{yes}{{\fontsize{20}{20}${\rm lifting}$}}
\psfrag{D}{{\fontsize{20}{20}$D$}}
\psfrag{p(D)}{{\fontsize{20}{20}$p(D)$}}
\psfrag{y}{{\fontsize{20}{20}$y$}}
\psfrag{v1}{{\fontsize{20}{20}$v_1$}}
\psfrag{v2}{{\fontsize{20}{20}$v_2$}}
\psfrag{v3}{{\fontsize{20}{20}$v_3$}}

\scalebox{0.45}{\includegraphics{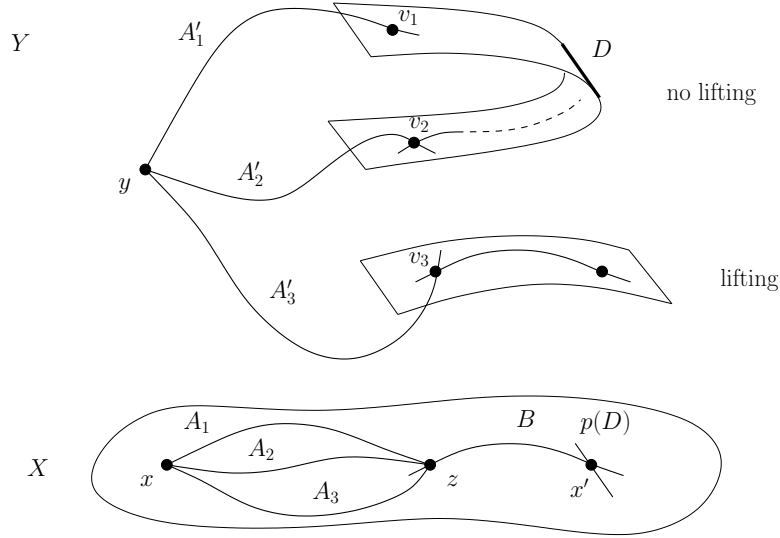}}
\caption{A non-liftable reducible curve}
\label{nonlift.pict}
\end{figure} 


There are only two points that need refinement
in this method. First, in positive characteristic
the local description of the ramification is much  more complicated than 
over $\c$. Second,
in order to make these choices more uniform, we should work with a
whole family of curves $B$. This actually also takes care of the
ramification problem. This is done next.

\end{say}

\begin{say}[Main construction]\label{main.fg.const}

Let $X$ be a smooth irreducible 
variety over an algebraically closed field $k$.
 Assume that we are given two families of
curves on $X$ 

$U\stackrel{u}{\leftarrow} A\stackrel{F}{\to} X$
and
$V\stackrel{v}{\leftarrow} B\stackrel{G}{\to} X$
with the following properties.
\begin{enumerate}
\item $U,V$ are irreducible $k$-varieties and
 $u,v$ are proper, smooth  morphisms with
irreducible  fibers.
\item $F$ is dominant.
\item $u$ has a section $s:U\to A$ such that  $F\circ s$ maps $U$
to a single point $x\in X(k)$.
\item  The images of the  maps
$\pi_1(A,s(u))\to \pi_1(X,x)$ for $u\in U$
topologically generate $\pi_1(X,x)$.
\item $G$ is smooth and for every divisor $D\subset X$ the general
$G(B_v)$ intersects $D$.
\end{enumerate}

Every point of the fiber product
$W:=A\times_XB$ can be thought of as a quadruplet
$(a_u\in A_u,b_v\in B_v)$ such that $F(a_u)=G(b_v)$.
To  this corresponds a reducible 
nodal curve $A_u\cup_{a_u\sim b_v}B_v$
obtained from the disjoint union of $A_u$ and  $B_v$ by identifying
the points $a_u$ and $b_v$.
Thus we have a flat family of reducible but connected
curves
$$
W{\leftarrow} C\stackrel{H}{\to} X.
$$
$s:U\to A$
gives a section $S:W\to C$ such that $H\circ S$ maps $W$ to $\{x\}$.
($C$ can be explicitly constructed as
$$
A\times_UA\times_XB \cup A\times_XB\times_VB
\subset A\times_UA\times_XB\times_VB,
$$
where the embedding is with $Id\times_X\Delta_B \cup \Delta_A\times_XId$
where $\Delta_B:B\to B\times_VB$ and $\Delta_A:A\to A\times_UA$
are the diagonal embeddings. The map $C\to W$ is projection to the
middle $A\times_XB$ of the above product. Finally $S$ is given by
the fiber product of $s$ with the identity map
$A\times_XB=U\times_UA\times_XB\to A\times_UA\times_XB$.)

\begin{prop} \label{2step.fg.prop}
Notation and assumptions as above.
 Then, for every open set $X^0\subset X$ containing $x$
$$
\pi_1(X^0,x)=
\langle \im[\pi_1(C_w\cap H^{-1}(X^0),s(w))\to \pi_1(X^0,x)]
: \forall w\in W(k) \rangle,
$$
where $\langle\cdots \rangle$ denotes the closed subgroup generated
by the images.
\end{prop}

Proof. It is enough to show  that if
 $(Y^0,y)\to (X^0,x)$ is a pointed finite \'etale cover
with $Y^0$ irreducible
and $p:Y\to X$ its extension to a (possibly ramified) finite morphism
then
$$
H_w:C_w=A_u\cup_{a_u\sim b_v}B_v\to X
$$
cannot be lifted  to $Y$ for some $w$.

If $Y\to X$ is unramified then there is a $u\in U(k)$ such that
$A_u\to X$ can not be lifted by assumption 4. Thus we are
left to consider the case when $Y\to X$ is ramified
and every $A_u\to X$ can  be lifted. We claim that in fact
 $F:A\to X$ lifts to $F_Y:A\to Y$ with $F_Y\circ s:U\to y$.
To see this, consider the fiber product $A\times_XY$.
$s:U\to A$ and the constant map $U\to y$ provide
$\sigma:U\to A\times_XY$. Let $A_Y\subset A\times_XY$  be the unique
irreducible component containing $\sigma(U)$. 
By assumption $A_Y\to A$ is an isomorphism on every fiber of $A\to U$,
hence an isomorphism. Thus $A_Y$ is the graph of the required lifting $F_Y$.

$F_Y$ is dominant, so there is an open set
$X^1\subset X^0$ such that  $F_Y(A)\supset  p^{-1}(X^1)$
and $X^1\subset G(B)$.

Let $z\in X^1(k)$ be a point and $p^{-1}(z)=\{v_1,\dots, v_m\}$.
By assumption there are $u_i\in U$ and $a_i\in A_{u_i}$
such that $F_Y(a_i)=v_i$. Pick any $v\in V(k)$ and $b\in B_v$
such that $G(b)=z$.

The quadruplets $(a_i\in A_{u_i}, b\in B_v)$ give
connected curves $C_i:=A_{u_i}\cup_{a_i\sim b}B_v$.
If $H_i:C_i\to X$ lifts to $H_{i,Y}:C_i\to Y$
then $H_{i,Y}(B_v)$ is  a lifting of $G:B_v\to X$
which passes through $v_i$. Thus we obtain:

\begin{claim} Under the above assumptions,  the fiber product
$B_v\times_XY$ has $m$ irreducible components, each
isomorphic to $B_v$.
\end{claim}

We  show that the this  is impossible, which implies that our assumption
is incorrect. Thus not all maps $H_i:C_i\to X$ lift to $H_{i,Y}:C_i\to Y$.
This will  complete the proof of (\ref{2step.fg.prop}).

It should be emphasized that at this point we only 
know that $B_v\times_XY$ has $m$ irreducible components, each
isomorphic to $B_v$. We do not know that $B_v\times_XY$
is a {\em disjoint} union of its irreducible components.
The key point is precisely to establish this for general $v\in V(k)$.

Set $B_Y:=B\times_XY$.
$G:B\to X$ is smooth, hence $B_Y\to Y$ is also smooth,
thus $B_Y$ is normal.

Consider the diagram
$$
\begin{array}{ccccc}
B_Y & \stackrel{p_B}{\to} & B& \stackrel{G}{\to} & X\\
v_Y\downarrow\hphantom{v_Y}  &&  \hphantom{v} \downarrow v&&\\
V & = & V.
\end{array}
$$

In characteristic zero, $B_Y$ normal implies that
the  geometric generic fiber
of $B_Y/V$ is also normal, hence it is a 
disjoint union of its irreducible components.

In positive characteristic, we have only a weaker result (\ref{npstein.lem}),
which still says that the  geometric generic fiber
of $B_Y/V$ is  a 
disjoint union of its irreducible components.
In particular,  the generic fiber of $B_Y/V$
is smooth over the generic fiber of $B/V$. Thus, after possibly
shrinking $V$, we may assume that $p_B$ is smooth.
$G$ is smooth by assumption, hence the composite
$B_Y\to X$ is also smooth.

 One can, however, factor this
as $B_Y\to Y\to X$. 
By the purity of branch loci, $Y\to X$ ramifies over a
whole divisor $D\subset X$,
 and so $B_Y\to X$ ramifies over every
point of $D$. This is only possible if the image of
$B_Y\to B\to X$ does not intersect $D$. 
This contradicts the assumption 5.
\qed
\end{say}

\begin{say}[Proof of (\ref{weakratlef})]

All that is left is  the construction of the families
$U\leftarrow A\to X$ and $V\leftarrow B\to X$ as in (\ref{main.fg.const}).
The rest is already taken care of by (\ref{2step.fg.prop}).

 For  $U\leftarrow A\to X$
we take  a family of rational curves through $x$
that covers $X$. Property 
(\ref{main.fg.const}.4) holds by (\ref{SRC.sc}).

In order to get the family $B\to V$, it
 is enough to get a very free curve $\p^1\to X$
which has positive intersection number with every divisor.
A free curve has nonnegative intersection number with every divisor.
Start with a family 
 $H:U\times \p^1\to X$ that
shows that $X$ is SRC. 
The only problem is that
its image may not intersect finitely many divisors
$D_1,\dots,D_s\subset X$. By \cite[29.5]{arko} 
there is a free morphism
 $g:\p^1\to X$ 
whose image intersects every $D_i$. 
Now take $g(\p^1)$ and a curve $H(\p^1_u)$ 
intersecting it.  By \cite[II.7]{Ko96}
their union can be smoothed to a  free
curve which has positive intersection number with every divisor.
The resulting morphism $G:B\to X$ is smooth by \cite[II.3.5.3]{Ko96}.
\qed
\end{say}

\begin{say}\label{main.thm.rem.pf}
 The additional properties listed
in (\ref{main.thm.rem}) are also easy to get.

First, 
by \cite[II.3.14]{Ko96}, $F$ is automatically an embedding
for general $w\in W$ if $\dim X\geq 3$, and we can
replace $W$ with an open subset at any time

Assume now that we want more ampleness from $F_w^*T_X$.
Take a family of morphisms $h:U\times \p^1\to X$
such that $h_u^*T_X(-d)$ is ample and $h$ is smooth.
Let us now look at all maps 
of one tooth combs
$g:C_0\cup C_1\to X$
where $g|_{C_0}$ is from the family $F$ and $g|_{C_1}$
from the family $h$. Smoothing these as in
(\ref{pf.of.23->4}), we get a new family $F_d:W_d\times \p^1\to X$
that also satisfies (\ref{main.thm.rem}.6)

Fix next an embedding $X\subset \p^N$. All smooth varieties
of given degree and dimension form a  bounded family
(the Chow variety).
Being SRC is an open property  (cf.\ \cite[IV.3.11]{Ko96}),
hence all SRC varieties 
of given degree and dimension form a  bounded family.
This implies that we can choose bounded  degree  curves
for the families $A,B$ in (\ref{main.fg.const}).
Thus we get families of bounded degree curves in (\ref{weakratlef}).

There is only one problem left, namely that during
the proof in (\ref{pf.of.23->4}) we have to take all
conjugates of a family of curves, and their number does depend
on the field. 
The families we have are open subsets of an irreducible component
of the Hom scheme $\Hom(\p^1,X)$.
We are dealing with maps of bounded degree
to a bounded family of varieties. The relative Hom scheme
of bounded degree maps is  quasiprojective
 (cf.\ \cite[I.1.10]{Ko96}),
 hence there is a  bound for the number of
irreducible families of bounded degree maps. The number of
conjugates cannot exceed this bound.\qed
\end{say}

\begin{say}[Proof of (\ref{ratlefx.cor})]

From (\ref{ratlefx.thm})
we know that $\pi_1(F^{-1}(X^0))\to G$ is surjective.
Applying (\ref{gen.semicont.prop}) to $F^{-1}(X^0)\to W$
we conclude that 
$\pi_1(\p^1_w\cap F^{-1}(X^0))\to G$ is surjective
for all $w$ in some open subset $W^G\subset W$.
\qed
\end{say}

\begin{say}[Proof of (\ref{ratlef.thm})]

Let us start with a variety $X$ and let $k(X)$ denote its field of
rational functions. The generic point  $x_g\in X$ is a $k(X)$-point of $X$,
thus we can apply (\ref{ratlefx.thm}) to $X_{k(X)}$
and  $x_g$. Thus we get
a geometrically irreducible $k(X)$-variety $V_g$ and a morphism
$$
F_g: V_g\times \p^1\to X_{k(X)}.
$$
$V_g$ can be thought of as the generic fiber of a map $V\to X$
where $V$ is a geometrically irreducible $k$-variety.
 $F_g$ extends to a map which
can be composed with the second  projection to obtain
$$
F_V:V\times \p^1\map X\times X\to X.
$$
$F_V$ need not be everywhere defined, but it
becomes a morphism after restriction to a suitable open set of the form
$V'\times \p^1$. By choosing $V'$ small enough we can also assume
that $F_V$ is free on $\{v\}\times \p^1$ for every $v\in V'$.

Let us now look at the relative Hom--scheme (cf.\ \cite[I.1.19]{Ko96})
$H:=\Hom_{V'}(V'\times \p^1, V'\times X)$. 
The restriction of $F_V$ to $V'$  determines a section
 $s:V'\to H$ and the universal morphism  $H\times\p^1\to V'\times X$
 is smooth over this
section by \cite[II.3.5.3]{Ko96}. This also imples that $H\to V'$ is
 smooth along $s(V')$, hence $H$ has a unique irreducible component
$W'\subset H$ which contains $s(V')$ and $W'$ is geometrically irreducible.
Furthermore, there is an open subset $W\subset W'$, containing $s(V')$
such that 
the universal morphism  $F:W\times\p^1\to V'\times X\to X$
is smooth. This gives the required morphism.
\qed

\end{say}

\begin{say}[Proof of (\ref{ratlef.cor2})]
\label{pf.of.rlc2}

Choose a partial compactification $Y\subset \bar Y$
such that $h$ extends to a proper morphism $\bar h:\bar Y\to X$.
By the upper semicontinuity of fiber dimensions,
there is a closed subset $T\subset X$ of codimension 2  
such that every fiber of $h$ over $X\setminus T$ has the same dimension.

This implies that  if $B$ is any irreducible curve and $B\to X\setminus T$
a morphism then every irreducible component of $B\times_XY$
dominates $B$. Indeed, the fiber product $B\times_XY$
is the preimage of the diagonal  under
the morphism $B\times Y\to X\times X$.
Therefore every irreducible component of $B\times_XY$ 
has dimension at least $\dim Y-\dim X+1$.
If $B$ maps to $X\setminus T$ then every fiber of
$B\times_XY\to B$ has dimension $\dim Y-\dim X$,
hence the claim.

Choose $F:W\times \p^1\to X$ as in (\ref{ratlef.thm}).
$F$ is smooth, so by an easy dimension count,
there is an open subset $W^1\subset W$ such that
$F(W^1\times \p^1)\subset X\setminus T$
(cf. \cite[9]{arko}).

Let
$$
Y^0:=p^{-1}(X^0)\stackrel{q}{\to}  Z^0\stackrel{r}{\to} X^0
$$
be as in (\ref{npstein.lem}), where $Z^0$ is normal and irreducible.
$q_w:
Y^0\times_X\p^1_w{\to}  Z^0\times_X\p^1_w $
 has geometrically irreducible fibers, thus it is enough to prove that
 $Z^0\times_X\p^1_w$ is irreducible  if $w\in W^0$.

$W^1\times\p^1\times_XZ^0$ is irreducible by
(\ref{ratlefx.thm}) and also smooth, and it has a section over $W^1$,
for instance  $w\mapsto (w,0,F(w,0))$.
So the generic fiber of $W^1\times\p^1\times_XZ^0\to W^1$
is irreducible over $k(W^1)$, smooth  and it also has a $k(W^1)$-point.
Therefore it is geometrically irreducible. 
Thus there is an open subset $W^0\subset W^1$ such that the fibers
of $W^1\times\p^1\times_XZ^0\to W^1$ are geometrically irreducible
 over $W^0$ (cf.\ \cite[IV.12]{ega}).
\qed
\end{say}

\section{Applications to non--closed fields}

In this section we derive some consequences of the previous results
to fields which are not algebraically closed.

\begin{rem}[Fundamental groups over non closed fields] 

If $k$ is a field then the connected \'etale covers of $\spec k$
are exactly the maps $\spec K\to \spec k$ where $K/k$
is a finite separable field extension. This implies that
$$
\pi_1(\spec (k))=\gal (k^{sep}/k)
$$
where $k^{sep}$ is the separable closure of $k$.
In general, the fundamental group of a $k$-scheme $X_k$  is 
made up of the  fundamental group of $X_{\bar k}$ and
the Galois group $\gal (k^{sep}/k)$.
To be precise, if $X_k$ is geometrically connected,
then there is an exact sequence
$$
1\to \pi_1(X_{\bar k},\bar x)\to  \pi_1(X_k,\bar x)\to \gal (k^{sep}/k)\to 1
$$
and every point of $X(k)$ defines a splitting of the sequence.
See \cite[IX.6.1]{sga1} for details.

This has an important consequence which indicates that
surjectivity assertions between fundamental groups
 are essentially geometric in nature:
\end{rem}

\begin{prop}\label{onto.is.geom}
 Let $k$ be a field and $f:Y\to X$  a morphism
of geometrically connected $k$-schemes. Let $\bar y\to Y$ be a 
geometric point.
Then the map of algebraic fundamental groups
$$
f_*: \pi_1(Y,\bar y)\to \pi_1(X,f(\bar y))
$$
is surjective (resp.
 its image has finite index)
iff the map of geometric fundamental groups
$$
f_*: \pi_1(Y_{\bar k},\bar y)\to \pi_1(X_{\bar k},f(\bar y))
$$
is surjective (resp.
 its image has finite index).\qed
\end{prop}

For arithmetic applications  of these results
it is useful to find maps  $\p^1\to X$
defined over the ground field. The current results
work very well for  certain fields:

\begin{say}[Large fields] We are interested in 
fields  $K$ which have the property that on any variety
with one smooth $K$-point the  $K$-points are
  Zariski dense.
Such fields  are called {\it large fields} in \cite{pop}.
The following
are some interesting classes of such fields:
\begin{enumerate}
\item Fields complete with respect to a discrete valuation.
This in particular includes the finite extensions of the $p$-adic
fields $\q_p$
and the power series fields $\f_q((t))$ over finite fields $\f_q$.
\item More generally, quotient fields of   local Henselian domains.
\item $\r$ and all real closed fields.
\item Infinite algebraic extensions of finite fields and, more
generally, pseudo algebraically closed fields, cf.\ \cite[Chap.\
10]{fj86}
\end{enumerate}
\end{say}

For large fields, the existence of a smooth $k$ point
in $\bar W$ in (\ref{ratlefx.thm})
implies that $W(k)$ is dense in $W$. 
Using (\ref{ratlefx.thm}) and (\ref{ratlefx.cor})
 we obtain the following:

\begin{thm} \label{ratlefx-k.thm}
Let $X$ be a smooth, projective, SRC variety over a large field $k$
and $x\in X(k)$.
 Then
 there is a dominant family of rational curves through $x$
(defined over $k$)
$$
F:W\times \p^1\to X,  \quad F(W\times (0{:}1))=\{x\}
$$
with the following properties:
\begin{enumerate}
\item $F_w^*T_X$ is ample for every $w\in W$.
\item Let $X^0\subset X$ be
any open $k$-subvariety containing $x$
and $\pi_1(X^0, x)\onto G$  a finite quotient
(whose kernel is open).
Then there is
a dense  subset of $k$-points
$W^0_G(k)\subset W$ such that for every $w\in W^0_G(k)$ the composed map
of algebraic fundamental groups
$$
\pi_1(\p_w^1\cap F^{-1}(X^0),(0{:}1))\to
\pi_1(X^0, x)\to G 
\qtq{is surjective.}\qed
$$ 
\end{enumerate}
\end{thm}

An arithmetic application  of these results is a
simple proof of the following theorem,
proved in increasing generality in
\cite{harb, ct-torsm, koll-rcfg}.
 \cite{mb01, mb02}  proved  analogous results
for finite and smooth group schemes in positive characteristic.

\begin{thm}\label{tors.gen.cor}
 Let $K$ be a characteristic zero large  field, 
$G$  a linear algebraic group scheme over $K$ and  $A$  a
principal homogeneous $G$-space.
 Then there is an open set
$V\subset \a^1_K$ containing $0$ and a geometrically irreducible
$G$-torsor $g:V_G\to V$ such that $g^{-1}(0)\cong A$
(as a $G$-space).
\end{thm}

Proof. Assume that $G$ acts on $A$ from the left
and choose an embedding $G\subset GL(n)$ over $K$. The group
 $A\times GL(n)$ admits a diagonal left action by $G$
and a right action by $GL(n)$ acting only on $GL(n)$. 

The right $G$-action makes the morphism
 $$
h:  Y^0:= G\backslash(A\times GL(n))\to 
 G\backslash(A\times GL(n))/G=:X^0
$$
 into a $G$-torsor.
Let $G^0\subset G$ be the connected component of the identity.
Then
 $$
Z^0:=G\backslash(A\times GL(n))/G^0\to X^0
$$
is a finite  \'etale cover, hence
it corresponds to a finite quotient
$\pi_1(X^0,x)\to H$ where 
 $x\in X(K)$ is the image of $G\backslash(A\times G)$. Let 
 $X\supset X^0$ be a smooth compactification of $X^0$. The variety
 $Z^0_{\bar K}$
is isomorphic to $GL(n)$, hence $X$ is unirational
and so SRC.
The fiber of $h$ over $x$ is isomorphic to
$A$. 

Apply (\ref{ratlefx-k.thm}) to $X$ with the point $x$
to get $f:(0\in\p^1)\to (x\in X)$ 
such that $\p^1\times_XZ^0$ is geometrically irreducible.
The fibers of $Y^0\to Z^0$ are torsors over $G^0$,
hence irreducible. Thus $\p^1\times_XY^0$ is also geometrically irreducible.
Set $V=f^{-1}(X^0)$ and $V_G=V\times_XY^0$.
\qed

\begin{rem}
We have  used the characteristic zero assumption in two places.
The first is the existence of a smooth compactification $X$.
Conjecturally this is always satisfied. A second point is that in general
we can conclude that $X$ is SRC only if $h:Y^0\to X^0$
is separable, that is when 
 $G$ is smooth. The computations of \cite{gille1, gille2}
suggest  that in general this may be  a quite subtle point.
\end{rem}

\begin{ack}   I   thank  E.\ Szab\'o
for helpful
comments and C.\ Araujo for the nice picture.
I am very grateful to  J.-L.\ Colliot-Th\'el\`ene for a  long list of
corrections and suggestions which improved the paper enormously.
Partial financial support was provided by  the NSF under grant number 
DMS-9970855 and DMS-02-00883. 
\end{ack}

\vskip1cm

\noindent Princeton University, Princeton NJ 08544-1000

\begin{verbatim}kollar@math.princeton.edu\end{verbatim}

\end{document}